\documentclass[12 pt]{article}%
\usepackage{amssymb,amsmath,amsfonts,amsthm}
\usepackage{amsmath}
\usepackage{amsfonts}
\usepackage{amssymb}
\usepackage{graphicx}%
\providecommand{\U}[1]{\protect\rule{.1in}{.1in}}
\theoremstyle{plain}
\newtheorem{thm}{Theorem}

\theoremstyle{definition}

\begin{document}
To appear in the Proceedings Volume of the International Symposium: Asymptotic
Methods in Stochastics, that was organized and held in honour of the work of
Miklos Csorgo on the occasion of his 80th birthday at Carleton University,
July 3-6, 2012. The Volume will be published in the Fields Institute's
Communication Series by Springer.

\bigskip

\begin{center}
{\LARGE Quenched Invariance Principles via Martingale Approximation}

\bigskip

Magda Peligrad\footnote{\textit{ }Supported in part by a Charles Phelps Taft
Memorial Fund grant, the NSA\ grant H98230-11-1-0135 and the NSF grant
DMS-1208237.}

\end{center}

Department of Mathematical Sciences, University of Cincinnati, PO Box 210025,
Cincinnati, Oh 45221-0025, USA. E-mail address: peligrm@ucmail.uc.edu

\bigskip

\begin{center}
\textbf{Abstract}

\bigskip

\end{center}

In this paper we survey the almost sure central limit theorem and its
functional form (quenched) for stationary and ergodic processes. For additive
functionals of a stationary and ergodic Markov chain these theorems are known
under the terminology of central limit theorem and its functional form,
started at a point. All these results have in common that they are obtained
via a martingale approximation in the almost sure sense. We point out several
applications of these results to classes of mixing sequences, shift processes,
reversible Markov chains, Metropolis Hastings algorithms.

\begin{center}
\textbf{ }
\end{center}

\section{Introduction and general considerations}

\quad In recent years there has been an intense effort towards a better
understanding of the structure and asymptotic behavior of stochastic
processes. For dependent sequences there are two basic techniques:
approximation with independent random variables or with martingales. Each of
these methods have its own strength. On one hand the processes that can be
treated by coupling with an independent sequence exhibit faster rates of
convergence in various limit theorems; on the other hand the class of
processes that can be treated by a martingale approximation is larger. There
are plenty of processes that benefit from approximation with a martingale.
Examples are: linear processes with martingale innovations, functions of
linear processes, reversible Markov chains, normal Markov chains, various
dynamical systems and the discrete Fourier transform of general stationary
sequences. A martingale approximation provides important information about
these structures because of their rich properties. They satisfy a broad range
of inequalities, they can be embedded into Brownian motion and they satisfy
various asymptotic results such as the functional conditional central limit
theorem and the law of the iterated logarithm. Moreover, martingale
approximation provides a simple and unified approach to asymptotic results for
many dependence structures. For all these reasons, in recent years martingale
approximation, "coupling with a martingale", has gained a prominent role in
analyzing dependent data. This is also due to important developments by
Liverani (1996), Maxwell-Woodroofe (2000), Derriennic-Lin (2001 a and b,
2003), Wu-Woodroofe (2004) and developments by Peligrad-Utev (2005),
Zhao-Woodroofe (2008 a and b), Voln\'{y} (2007), Peligrad-Wu (2010) among
others. Many of these new results, originally designed for Markov operators,
(see Kipnis-Varadhan, 1986; Derriennic-Lin, 2001 b) have made their way into
limit theorems for stochastic processes.

This method has been shown to be well suited to transport from the martingale
to the stationary process either the conditional central limit theorem or
conditional invariance principle in probability. As a matter of fact, papers
by Dedecker-Merlev\`{e}de-Voln\'{y} (2007), Zhao and Woodroofe (2008 b),
Gordin and Peligrad (2011), point out characterizations of stochastic
processes that can be approximated by martingales in quadratic mean. These
results are useful to treat evolutions in "annealed" media.

In this survey we address the question of limit theorems started at a point
for almost all points. These type of results are also known under the name of
quenched limit theorems or almost sure conditional invariance principles.
Limit theorems for stochastic processes that do not start from equilibrium is
timely and motivated by recent development in evolutions in quenched random
environment, random walks in random media, for instance as in Rassoul-Agha and
Sepp\"{a}l\"{a}inen (2007). Moreover recent discoveries by Voln\'{y} and
Woodroofe (2010 a) show that many of the central limit theorems satisfied by
classes of stochastic processes in equilibrium, fail to hold when the
processes are started from a point. Special attention will be devoted to
normal and reversible Markov chains and several results and open problems will
be pointed out. These results are very important since reversible Markov
chains have applications to statistical mechanics and to Metropolis Hastings
algorithms used in Monte Carlo simulations. The method of proof of this type
of limiting results are approximations with martingale in an almost sure sense.

The field of limit theorems for stochastic processes is closely related to
ergodic theory and dynamical systems. All the results for stationary sequences
can be translated in the language of Markov operators.

\section{Limit theorems started at a point via martingale approximation}

\quad In this section we shall use the framework of strictly stationary
sequences adapted to a stationary filtrations that can be introduced in
several equivalent ways, either by using a measure preserving transformation
or as a functional of a Markov chain with a general state space. It is just a
difference of language to present the theory in terms of stationary processes
or functionals of Markov chains.

Let $(\Omega,\mathcal{A},{\mathbb{P}})$ be a probability space, and
$T:\Omega\mapsto\Omega$ be a bijective bimeasurable transformation preserving
the probability ${\mathbb{P}}$. A set $A\in\mathcal{A}$ is said to be
invariant if $T(A)=A$. We denote by $\mathcal{I}$ the $\sigma$-algebra of all
invariant sets. The transformation $T$ is ergodic with respect to probability
${\mathbb{P}}$ if each element of $\mathcal{I}$ has measure $0$ or $1$. Let
$\mathcal{F}_{0}$ be a $\sigma$-algebra of $\mathcal{A}$ satisfying
$\mathcal{F}_{0}\subseteq T^{-1}(\mathcal{F}_{0})$ and define the
nondecreasing filtration $(\mathcal{F}_{i})_{i\in{\mathbb{Z}}}$ by
$\mathcal{F}_{i}=T^{-i}(\mathcal{F}_{0})$. Let $X_{0}$ be a $\mathcal{F}_{0}%
$-measurable, square integrable and centered random variable. Define the
sequence $(X_{i})_{i\in\mathbb{Z}}$ by $X_{i}=X_{0}\circ T^{i}$. Let
$S_{n}=X_{1}+\cdots+X_{n}.$ For $p\geq1,\Vert.\Vert_{p}$ denotes the norm in
$\mathbb{L}_{p}(\Omega,\mathcal{A},{\mathbb{P}}).$ In the sequel we shall
denote by $\mathbb{E}_{0}(X)=\mathbb{E}(X|\mathcal{F}_{0})$.

\quad The conditional central limit theorem plays an essential role in
probability theory and statistics. It asserts that the central limit theorem
holds in probability under the measure conditioned by the past of the process.
More precisely this means that for any function $f$ which is continuous and
bounded we have%

\begin{equation}
\mathbb{E}_{0}(f(S_{n}/\sqrt{n}))\rightarrow\mathbb{E}(f(\sigma N))\text{ in
probability,} \label{CCLT}%
\end{equation}
where $N$ is a standard normal variable and $\sigma$ is a positive constant.
Usually we shall have the interpretation $\sigma^{2}=\lim_{n\rightarrow\infty
}$\textrm{var}$(S_{n})/n.$

This conditional form of the CLT is a stable type of convergence that makes
possible the change of measure with a majorizing measure, as discussed in
Billingsley (1968), Rootz\'{e}n (1976), and Hall and Heyde (1980).
Furthermore, if we consider the associated stochastic process
\[
W_{n}(t)=\frac{1}{\sqrt{n}}S_{[nt]}\text{,}%
\]
where $[x]$ denotes the integer part of $x$, then the conditional CLT\ implies
the convergence of the finite dimensional distributions of $W_{n}(t)$ to those
of $\sigma W(t)$ where $W(t)$ is the standard Brownian Motion; this
constitutes an important step in establishing the functional CLT (FCLT). Note
that $W_{n}(t)$ belongs to the space $D[0,1]$, the set of functions on $[0,1]$
which are right continuous and have left hands limits. We endow this space
with the uniform topology.

By the conditional functional central limit theorem we understand that for any
function $f$ continuous and bounded on $D[0,1]$ we have
\begin{equation}
\mathbb{E}_{0}(f(W_{n}))\rightarrow\mathbb{E}(f(\sigma W))\ \text{in
probability.} \label{CIP}%
\end{equation}
There is a considerable amount of research concerning this problem. We mention
papers by Dedecker and Merlev\`{e}de (2002), Wu and Woodroofe (2004) and Zhao
and Woodroofe (2008 b) among others.

The quenched versions of these theorems are obtained by replacing the
convergence in probability by convergence almost sure. In other words the
almost sure conditional theorem states that, on a set of probability one, for
any function $f$ which is continuous and bounded we have%
\begin{equation}
\mathbb{E}_{0}(f(S_{n}/\sqrt{n}))\rightarrow\mathbb{E}(f(\sigma N))\text{,}
\label{asCLT}%
\end{equation}
while by almost sure conditional functional central limit theorem we
understand that, on a set of probability one, for any function $f$ continuous
and bounded on $D[0,1]$ we have
\begin{equation}
\mathbb{E}_{0}(f(W_{n}))\rightarrow\mathbb{E}(f(\sigma W))\text{.}
\label{asCIP}%
\end{equation}

We introduce now the stationary process as a functional of a Markov chain.

We assume that $(\xi_{n})_{n\in\mathbb{Z}}$ is a stationary ergodic Markov
chain defined on a probability space $(\Omega,\mathcal{F},\mathbb{P})$ with
values in a Polish space $(S,\mathcal{S})$. The marginal distribution is
denoted by $\pi(A)=\mathbb{P}(\xi_{0}\in A),$ $A\in\mathcal{S}.$ Next, let
$\mathbb{L}_{2}^{0}(\pi)$ be the set of functions $h$ such that $||h||_{2,\pi
}^{2}=\int h^{2}d\pi<\infty$ and $\int hd\pi=0.$ Denote by $\mathcal{F}_{k}$
the $\sigma$--field generated by $\xi_{j}$ with $j\leq k$, $X_{j}=h(\xi_{j})$.
Notice that any stationary sequence $(Y_{k})_{k\in\mathbb{Z}}$ can be viewed
as a function of a Markov process $\xi_{k}=(Y_{j};j\leq k)$ with the function
$g(\xi_{k})=Y_{k}$. Therefore the theory of stationary processes can be
imbedded in the theory of Markov chains.

In this context by the central limit theorem started at a point (quenched) we
understand the following fact: let $\mathbb{P}^{x}$ be the probability
associated with the process started from $x$ and let $\mathbb{E}^{x}$ be the
corresponding expectation. Then, for $\pi-$almost every $x$, for every
continuous and bounded function $f$,
\begin{equation}
\mathbb{E}^{x}(f(S_{n}/\sqrt{n}))\rightarrow\mathbb{E}(f(\sigma N))\text{.}
\label{PxCLT}%
\end{equation}
By the functional CLT started at a point we understand that, for $\pi-$almost
every $x$, for every function $f$ continuous and bounded on $D[0,1]$,
\begin{equation}
\mathbb{E}^{x}(f(W_{n}))\rightarrow\mathbb{E}(f(\sigma W))\text{.}
\label{PxIP}%
\end{equation}
where, as before $W$ is the standard Brownian motion on $[0,1].$

It is remarkable that a martingale with square integrable stationary and
ergodic differences satisfies the quenched CLT in its functional form. For a
complete and careful proof of this last fact we direct to Derriennic and Lin
(2001 a, page 520). This is the reason why a fruitful approach to find classes
of processes for which quenched limit theorems hold is to approximate partial
sums by a martingale.

The martingale approximation as a tool in studying the asymptotic behavior of
the partial sums $S_{n}$ of stationary stochastic processes goes back to
Gordin (1969) who proposed decomposing the original stationary sequence into a
square integrable stationary martingale $M_{n}=\sum_{i=1}^{n}D_{i}$ adapted to
$(\mathcal{F}_{n})$, such that $S_{n}=M_{n}+R_{n}$ where $R_{n}$ is a
telescoping sum of random variables, with the basic property that $\sup
_{n}||R_{n}||_{2}<\infty$.

For proving conditional CLT for stationary sequences, a weaker form of
martingale approximation was pointed out by many authors (see for instance
Merlev\`{e}de-Peligrad-Utev, 2006 for a survey).

An important step forward was the result by Heyde (1974) who found sufficient
conditions for the decomposition
\begin{equation}
S_{n}=M_{n}+R_{n}\text{ with }\mathbb{\ }R_{n}/\sqrt{n}\rightarrow0\text{ in
}\mathbb{L}_{2}\text{.} \label{approxl2}%
\end{equation}
Recently, papers by Dedecker-Merlev\`{e}de-Voln\'{y} (2007) and by
Zhao-Woodroofe (2008 b) deal with necessary and sufficient conditions for
martingale approximation with an error term as in (\ref{approxl2}).

The approximation of type (\ref{approxl2}) is important since it makes
possible to transfer from martingale the conditional CLT defined in
(\ref{CCLT}), where $\sigma=||D_{0}||_{2}.$

The theory was extended recently in Gordin-Peligrad (2011) who developed
necessary and sufficient conditions for a martingale decomposition with the
error term satisfying
\begin{equation}
\max_{1\leq j\leq n}|S_{j}-M_{j}|/\sqrt{n}\rightarrow0\text{ in }%
\mathbb{L}_{2}\text{.} \label{maxapprox}%
\end{equation}
This approximation makes possible the transport from the martingale to the
stationary process the conditional functional central limit theorem stated in
(\ref{CIP}). These results were surveyed in Peligrad (2010).

The martingale approximation of the type (\ref{maxapprox}) brings together
many disparate examples in probability theory. For instance, it is satisfied
under Hannan (1973, 1979) and Heyde (1974) projective condition.
\begin{equation}
\mathbb{E}(X_{0}|\mathcal{F}_{-\infty})=0\quad\text{almost surely and}%
\quad\sum_{i=1}^{\infty}\Vert\mathbb{E}_{-i}(X_{0})-\mathbb{E}_{-i-1}%
(X_{0})\Vert_{2}<\infty\text{;} \label{PROJ}%
\end{equation}
It is also satisfied for classes of mixing processes; additive functionals of
Markov chains with normal or symmetric Markov operators.

A very important question is to establish quenched version of conditional CLT
and conditional FCLT, i.e. the invariance principles as in (\ref{asCLT}) and
also in (\ref{asCIP}) (or equivalently as in (\ref{PxCLT}) and also in
(\ref{PxIP})). There are many examples of stochastic processes satisfying
(\ref{maxapprox}) for which the conditional CLT does not hold in the almost
sure sense. For instance condition (\ref{PROJ}) is not sufficient for
(\ref{asCLT}) as pointed out by Voln\'{y} and Woodroofe (2010 a). In order to
transport from the martingale to the stationary process the almost sure
invariance principles the task is to investigate the approximations of types
(\ref{approxl2}) or (\ref{maxapprox}) with an error term well adjusted to
handle this type of transport. These approximations should be of the type, for
every $\varepsilon>0$
\begin{equation}
\mathbb{P}_{0}[|S_{n}-M_{n}|/\sqrt{n}>\varepsilon]\rightarrow0\text{
}a.s.\text{ or }\mathbb{\ P}_{0}[\max_{1\leq i\leq n}|S_{i}-M_{i}|/\sqrt
{n}>\varepsilon]\rightarrow0\text{ }a.s. \label{a.s. aprox}%
\end{equation}
where $(M_{n})_{n}$ is a martingale with stationary and ergodic differences
and we used the notation $\mathbb{P}_{0}(A)=\mathbb{P}(A|\mathcal{F}_{0})$.
They are implied in particular by stronger approximations such that%
\[
|S_{n}-M_{n}|/\sqrt{n}\rightarrow0\text{ }a.s.\text{ or }\max_{1\leq i\leq
n}|S_{i}-M_{i}|/\sqrt{n}\rightarrow0\text{ }a.s.
\]
Approximations of these types have been considered in papers by Zhao-Woodroofe
(2008 a), Cuny (2011), Merlev\`{e}de- Peligrad M.- Peligrad C. (2011) among others.

In the next subsection we survey resent results and point out several classes
of stochastic processes for which approximations of the type (\ref{a.s. aprox}%
) hold.

For cases where a stationary martingale approximation does not exist or cannot
be pointed out, a nonstationary martingale approximation is a powerful tool.
This method was occasionally used to analyze a stochastic process. Many ideas
are helpful in this situation ranging from simple projective decomposition of
sums as in Gordin and Lifshitz (1981) to more sophisticated tools. One idea is
to divide the variables into blocks and then to approximate the sums of
variables in each block by a martingale difference, usually introducing a new
parameter, the block size, and changing the filtration. This method was
successfully used in the literature by Philipp-Stout (1975), Shao (1995),
Merlev\`{e}de-Peligrad (2006), among others. Alternatively, one can proceed as
in Wu-Woodroofe (2004), who constructed a nonstationary martingale
approximation for a class of stationary processes without partitioning the
variables into blocks.

Recently Dedecker-Merlev\`{e}de-Peligrad (2012) used a combination of blocking
technique and a row-wise stationary martingale decomposition in order to
enlarge the class of random variables known to satisfy the quenched invariance
principles. To describe this approach, roughly speaking, one considers an
integer $m=m(n)$ large but such that $n/m\rightarrow\infty.$ Then one forms
the partial sums in consecutive blocks of size $m$, $Y_{j}^{n}=X_{m(j-1)+1}%
+....+X_{mj},\ 1\leq j\leq k,$ $k=[n/m].$ Finally, one considers the
decomposition%
\begin{equation}
S_{n}=M_{n}^{n}+R_{n}^{n}, \label{LBM}%
\end{equation}
where $M_{n}^{n}=%
{\displaystyle\sum\limits_{j=1}^{n}}
D_{j}^{n},$ with $D_{j}^{n}=Y_{j}^{n}-E(Y_{j}^{n}|\mathcal{F}_{m(j-1)})$ a
triangular array of row-wise stationary martingale differences.

\subsection{ Functional Central limit theorem started at a point under
projective criteria.}

\quad We have commented that condition (\ref{PROJ}) is not sufficient for the
validity of the almost sure CLT started from a point. Here is a short history
of the quenched CLT under projective criteria.\textbf{\ }A result in Borodin
and Ibragimov (1994, ch.4, section \ 8) states that if $||\mathbb{E}_{0}%
(S_{n})||_{2}$ is bounded, then the CLT in its functional form started at a
point (\ref{asCIP}) holds. Later, Derriennic-Lin (2001 a and b, 2003) improved
on this result imposing the condition $||\mathbb{E}_{0}(S_{n})||_{2}%
=O(n^{1/2-\epsilon})$ with $\epsilon>0$ (see also Rassoul-Agha and
Sepp\"{a}l\"{a}inen, 2008). A step forward was made by Cuny (2011) who
improved the condition to $||\mathbb{E}_{0}(S_{n})||_{2}=O(n^{1/2}(\log
n)^{-2}(\log\log n)^{-1-\delta})$ with $\delta>0$, by using sharp results on
ergodic transforms in Gaposhkin (1996).

We shall describe now the recent progress made on the functional central limit
theorem started at a point under projective criteria. We give here below three
classes of stationary sequences of centered square integrable random variables
for which both quenched central limit theorem and its quenched functional form
given in (\ref{asCLT}) and (\ref{asCIP}) hold with $\sigma^{2}=\lim
_{n\rightarrow\infty}$\textrm{var}$(S_{n})/n$, provided the sequences are
ergodic. If the sequences are not ergodic then then the results still hold but
with $\sigma^{2}$ replaced by the random variable $\eta$ described as
$\eta=\lim_{n\rightarrow\infty}\mathbb{E}(S_{n}^{2}|\mathcal{I})/n$ and
$\mathbb{E(\eta)=}\sigma^{2}.$ For simplicity we shall formulate the results
below only for ergodic sequences.

1. \textbf{Hannan}-\textbf{Heyde projective criterion}. Cuny-Peligrad (2012)
(see also Voln\'{y}-Woodroofe, 2010 b) showed that (\ref{asCLT}) holds under
the condition
\begin{equation}
\frac{\mathbb{E}(S_{n}|\mathcal{F}_{0})}{\sqrt{n}}\rightarrow0\quad
\text{almost surely and}\quad\sum_{i=1}^{\infty}\Vert\mathbb{E}_{-i}%
(X_{0})-\mathbb{E}_{-i-1}(X_{0})\Vert_{2}<\infty\text{.} \label{CuPe}%
\end{equation}
The functional form of this result was established in Cuny-Voln\'{y} (2013).

2. \textbf{Maxwell and Woodroofe condition}. The convergence in (\ref{asCIP})
holds under Maxwell-Woodroofe (2000) condition,
\begin{equation}
\sum_{k=1}^{\infty}\frac{||\mathbb{E}_{0}(S_{k})||_{2}}{k^{3/2}}<\infty,
\label{MW}%
\end{equation}
as recently shown in Cuny-Merlev\`{e}de (2012). In particular both conditions
(\ref{CuPe}) and (\ref{MW}) and is satisfied if
\begin{equation}
\sum_{k=1}^{\infty}\frac{||\mathbb{E}_{0}(X_{k})||_{2}}{k^{1/2}}<\infty.
\label{MWcor}%
\end{equation}

3. \textbf{Dedecker-Rio condition.} In a recent paper
Dedecker-Merlev\`{e}de-Peligrad (2012) proved (\ref{asCIP}) under the
condition%
\begin{equation}
\sum_{k\geq0}\Vert X_{0}{\mathbb{E}}_{0}(X_{k})\Vert_{1}<\infty. \label{Rio}%
\end{equation}

The first two results were proved using almost sure martingale approximation
of type (\ref{a.s. aprox}). The third one was obtained using the large block
method described in (\ref{LBM}).

Papers by Durieu-Voln\'{y} (2008) and Durieu (2009) suggest that conditions
(\ref{CuPe}), (\ref{MW}) and (\ref{Rio}) are independent. They have different
areas of applications and they lead to optimal results in all these
applications. Condition (\ref{CuPe}) is well adjusted for linear processes. It
was shown in Peligrad and Utev (2005) that the Maxwell-Woodroofe condition
(\ref{MW}) is satisfied by $\rho-$mixing sequences with logarithmic rate of
convergence to $0$. Dedecker-Rio (2000) have shown that condition (\ref{Rio})
is verified for strongly mixing processes under a certain condition combining
the tail probabilities of the individual summands with the size of the mixing
coefficients. For example, one needs a polynomial rate on the strong mixing
coefficients when moments higher than two are available. However, the classes
described by projection conditions have a much larger area of applications
than mixing sequences. They can be verified by linear processes and dynamical
systems that satisfy only weak mixing conditions (Dedecker-Prieur 2004 and
2005, Dedecker-Merlev\`{e}de-Peligrad (2012) among others). More details about
the applications are given in Section 3.

Certainly, these projective conditions can easily be formulated in the
language of Markov operators by using the fact that ${\mathbb{E}}_{0}%
(X_{k})=Q(f)(\xi_{0}).$ In this language $\mathbb{E}_{0}(S_{k})=(Q+Q^{2}%
+...+Q^{k})(f)(\xi_{0}).$

\subsection{Functional Central limit theorem started at a point for normal and
reversible Markov chains.}

\quad In 1986 Kipnis and Varadhan proved the functional form of the central
limit theorem as in (\ref{CIP}) for square integrable mean zero additive
functionals $f\in\mathbb{L}_{2}^{0}(\pi)$ of stationary reversible ergodic
Markov chains $(\xi_{n})_{n\in\mathbb{Z}}$ with transition function $Q(\xi
_{0},A)=P(\xi_{1}\in A|\xi_{0})$ under the natural assumption $var(S_{n})/n$
is convergent to a positive constant. This condition has a simple formulation
in terms of spectral measure $\rho_{f}$ of the function $f$ with respect to
self-adjoint operator $Q$ associated to the reversible Markov chain, namely
\begin{equation}
\int\nolimits_{-1}^{1}\frac{1}{1-t}\rho_{f}(dt)<\infty\text{.} \label{SR1}%
\end{equation}
This result was established with respect to the stationary probability law of
the chain. (Self-adjoint means $Q=Q^{\ast},$ where $Q$ also denotes the
operator $Qf(\xi)=\int f(x)Q(\xi,dx)$; $Q^{\ast}$ is the adjoint operator
defined by $<Qf,g>=<f,Q^{\ast}g>$,$\ $for{\ every }${f}${\ and }${g}${\ in
$\mathbb{L}_{2}(\pi)$}).

The central limit theorem (\ref{CCLT}) for stationary and ergodic Markov
chains with normal operator $Q$ ($QQ^{\ast}=Q^{\ast}Q)$, holds under a similar
spectral assumption, as discovered by Gordin-Lifshitz (1981) (see also and
Borodin-Ibragimov, 1994, ch. 4 sections 7-8). A sharp sufficient condition in
this case in terms of spectral measure is%

\begin{equation}
\int\nolimits_{D}\frac{1}{|1-z|}\rho_{f}(dz)<\infty\text{.} \label{SN}%
\end{equation}
where $D$ is the unit disk.

Examples of reversible Markov chains frequently appear in the study of
infinite systems of particles, random walks or processes in random media. A
simple example of a normal Markov chain is a random walk on a compact group.
Other important example of reversible Markov chain is the extremely versatile
(independent) Metropolis Hastings Algorithm which is the modern base of Monte
Carlo simulations.

An important problem is to investigate the validity of the almost sure central
limit theorem started at a point for stationary ergodic normal or reversible
Markov chains. As a matter of fact, in their remark (1.7), Kipnis-Varadhan
(1986) raised the question if their result also holds with respect to the law
of the Markov chain started from $x$, for almost all $x,$ as in (\ref{PxIP}).

\textbf{Conjecture: }For any square integrable mean $0$ function of reversible
Markov chains satisfying condition (\ref{SR1}) the functional central limit
theorem started from a point holds for almost all points. The same question is
raised for continuous time reversible Markov chains.

The answer to this question for reversible Markov chains with continuous state
space is still unknown and has generated a large amount of research. The
problem of quenched CLT for normal stationary and ergodic Markov chains was
considered by Derriennic-Lin (2001 a) and Cuny (2011), among others, under
some reinforced assumptions on the spectral condition. Concerning normal
Markov chains, Derriennic-Lin (2001 a) pointed out that the central limit
theorem started at a point does not hold for almost all points under condition
(\ref{SN}). Furthermore, Cuny-Peligrad (2012) proved that there is a
stationary and ergodic normal Markov chain and a function $f\in\mathbb{L}%
_{2}^{0}(\pi)$ such that
\[
\int\nolimits_{D}\frac{|\log(|1-z|)\log\log(|1-z|)|}{|1-z|}\rho_{f}(dz)<\infty
\]
and such that the central limit theorem started at a point fails, for $\pi
-$almost all starting points.

However the condition
\begin{equation}
\int\nolimits_{-1}^{1}\frac{(\log^{+}|\log(1-t)|)^{2}}{1-t}\rho_{f}%
(dt)<\infty, \label{SN2}%
\end{equation}
is sufficient to imply central limit theorem started at a point (\ref{PxCLT})
for reversible Markov chains for $\pi-$almost all starting points. Note that
this condition is a slight reinforcement of condition (\ref{SN}).

It is interesting to note that by Cuny (2011, Lemma 2.1), condition
(\ref{SN2}) is equivalent to the following projective criterion%
\begin{equation}
\sum\nolimits_{n}\frac{(\log\log n)^{2}||\mathbb{E}_{0}(S_{n})||_{2}^{2}%
}{n^{2}}<\infty\text{.} \label{log}%
\end{equation}
Similarly, condition (\ref{SN}) in the case where $Q$ is symmetric, is
equivalent to
\begin{equation}
\sum\nolimits_{n}\frac{||\mathbb{E}_{0}(S_{n})||_{2}^{2}}{n^{2}}%
<\infty\text{.} \label{CLT-rev}%
\end{equation}

\section{Applications}

\quad Here we list several classes of stochastic processes satisfying quenched
CLT and quenched invariance principles. They are applications of the results
given in Section 2.

\subsection{Mixing processes}

\quad In this subsection we discuss two classes of mixing sequences which are
extremely relevant in the study of Markov chains, Gaussian processes and
dynamical systems.

We shall introduce the following mixing coefficients: For any two $\sigma
$-algebras $\mathcal{A}$ and $\mathcal{B}$ define the strong mixing
coefficient $\alpha(\mathcal{A}$,$\mathcal{B)}$:%

\[
\alpha(\mathcal{A},\mathcal{B)=}\sup\{|\mathbb{P}(A\cap B)-\mathbb{P}%
(A)\mathbb{P}(B)|;A\in\mathcal{A},B\in\mathcal{B\}}.
\]
The $\rho-$mixing coefficient, known also under the name of maximal
coefficient of correlation $\rho(\mathcal{A}$,$\mathcal{B})$ is defined as:
\[
\rho(\mathcal{A},\mathcal{B})=\sup\{\mathrm{Cov}(X,Y)/\Vert X\Vert_{2}\Vert
Y\Vert_{2}\;:\;X\in\mathbb{L}_{2}(\mathcal{A}),\text{ }Y\in\mathbb{L}%
_{2}(\mathcal{B})\}.
\]
For the stationary sequence of random variables $(X_{k})_{k\in\mathbb{Z}},$ we
also define $\mathcal{F}_{m}^{n}$ the $\sigma$--field generated by $X_{i}$
with indices $m\leq i\leq n,$ $\mathcal{F}^{n}$ denotes the $\sigma$--field
generated by $X_{i}$ with indices $i\geq n,$ and $\mathcal{F}_{m}$ denotes the
$\sigma$--field generated by $X_{i}$ with indices $i\leq m.$ The sequences of
coefficients $\alpha(n)$ and $\rho(n)$ are then defined by
\[
\alpha(n)=\alpha(\mathcal{F}_{0},\mathcal{F}^{n}\mathcal{)},\ \rho
(n)=\rho(\mathcal{F}_{0},\mathcal{F}^{n}\mathcal{)}.
\]
\quad An equivalent definition for $\rho(n)$ is
\begin{equation}
\rho(n)=\sup\{\Vert\mathbb{E}(Y|\mathcal{F}_{0})\Vert_{2}/\Vert Y\Vert
_{2}:\;Y\in\mathbb{L}_{2}(\mathcal{F}^{n}),\mathbb{E}(Y)=0\}. \label{rorick}%
\end{equation}
\quad Finally we say that the stationary sequence is strongly mixing if
$\alpha(n)\rightarrow0$ as $n\rightarrow\infty$, and $\rho-$mixing if
$\rho(n)\rightarrow0$ as $n\rightarrow\infty$. It should be mentioned that a
$\rho-$mixing sequence is strongly mixing. Furthermore, a stationary strongly
mixing sequence is ergodic. For an introduction to the theory of mixing
sequences we direct the reader to the books by Bradley (2007).

In some situations weaker forms of strong and $\rho-$mixing coefficients can
be useful, when $\mathcal{F}^{n}$ is replaced by the sigma algebra generated
by only one variable, $X_{n},$ denoted by $\mathcal{F}_{n}^{n}.$ We shall use
the notations $\tilde{\alpha}(n)=\alpha(\mathcal{F}_{0}$,$\mathcal{F}_{n}%
^{n}\mathcal{)}$ and $\tilde{\rho}(n)=\rho(\mathcal{F}_{0}$,$\mathcal{F}%
_{n}^{n}\mathcal{)}$.

By verifying the conditions in Section 3, we can formulate:

\begin{thm}
\label{mixing}Let $(X_{n})_{n\in\mathbb{Z}}$ be a stationary and ergodic
sequence of centered square integrable random variables. The quenched CLT\ and
its quenched functional form as in (\ref{asCLT}) and (\ref{asCIP}) hold with
$\sigma^{2}=\lim_{n\rightarrow\infty}$\textrm{var}$(S_{n})/n$ under one of the
following three conditions:%
\begin{equation}
\sum_{k=1}^{\infty}\frac{\tilde{\rho}(k)}{\sqrt{k}}<\infty.
\label{cond ro star}%
\end{equation}%
\begin{equation}
\sum_{k=1}^{\infty}\frac{\rho(k)}{k}<\infty. \label{condrho}%
\end{equation}%
\begin{equation}
\sum_{k=1}^{\infty}\int_{0}^{\tilde{\alpha}(k)}Q^{2}(u)du<\infty,
\label{cond alpha}%
\end{equation}
where $Q$ denotes the generalized inverse of the function $t\rightarrow
\mathbb{P}(|X_{0}|>t)$.
\end{thm}

We mention that under condition (\ref{condrho}) the condition of ergodicity is
redundant. Also if (\ref{cond alpha}) holds with $\tilde{\alpha}(k)$ replaced
by $\alpha(k),$ then the sequence is again ergodic.

In order to prove this theorem under (\ref{cond ro star}) one verifies
condition (\ref{MWcor}) via the estimate%
\[
\mathbb{E}(\mathbb{E}_{0}(X_{k}))^{2}=\mathbb{E}(X_{k}\mathbb{E}_{0}%
(X_{k}))\leq\tilde{\rho}(k)||X_{0}||_{2}^{2},
\]
which follows easily from the definition of $\tilde{\rho}.$

Condition (\ref{condrho}) is used to verify condition (\ref{MW}). This was
verified in the Peligrad-Utev-Wu (2007) via the inequalities%
\[
\Vert\mathbb{E}(S_{2^{r+1}}|\mathcal{F}_{0})\Vert_{2}\leq c\sum_{j=0}%
^{r}2^{j/2}\rho(2^{j})
\]
and
\begin{equation}
\sum_{r=0}^{\infty}{\frac{\Vert\mathbb{E}(S_{2^{r}}|\mathcal{F}_{0})\Vert_{2}%
}{2^{r/2}}}\leq c\sum_{j=0}^{\infty}\rho(2^{j})<\infty\text{.} \label{diatic}%
\end{equation}
Furthermore (\ref{diatic}) easily implies (\ref{MW}). For more details on this
computation we also direct the reader to the survey paper by
Merlev\`{e}de-Peligrad-Utev (2006).

To get the quenched results under condition (\ref{cond alpha}) the condition
(\ref{Rio}) is verified via the following identity taken from Dedecker-Rio
(2000, (6.1))
\begin{equation}
\mathbb{E}|X_{0}\mathbb{E}(X_{k}|\mathcal{F}_{0})|=\mathrm{Cov}(|X_{0}%
|(I_{\{\mathbb{E}(X_{k}|\mathcal{F}_{0})>0\}}-I_{\{\mathbb{E}(X_{k}%
|\mathcal{F}_{0})\leq0\}}),X_{k}). \label{cov}%
\end{equation}
By applying now Rio's (1993) covariance inequality we obtain%
\[
\mathbb{E}|X_{0}\mathbb{E}(X_{k}|\mathcal{F}_{0})|\leq c\int_{0}%
^{\tilde{\alpha}(k)}Q^{2}(u)du.
\]

It is obvious that condition (\ref{cond ro star}) requires a polynomial rate
of convergence to $0$ of $\tilde{\rho}(k);$ condition (\ref{condrho}) requires
only a logarithmic rate for $\rho(n).$ To comment about condition
(\ref{cond alpha}) it is usually used in the following two forms:

-either the variables are almost sure bounded by a constant, and then the
requirement is $\sum_{k=1}^{\infty}\tilde{\alpha}(k)<\infty.$

-the variables have finite moments of order $2+\delta$ for some $\delta>0,$
and then the condition on mixing coefficients is $\sum_{k=1}^{\infty
}k^{2/\delta}\tilde{\alpha}(k)<\infty.$

\subsection{Shift processes.}

\quad In this sub-section we apply condition (\ref{MW}) to linear processes
which are not mixing in the sense of previous subsection. This class is known
under the name of one-sided shift processes, also known under the name of
Raikov sums.

Let us consider a Bernoulli shift. Let $\{{\varepsilon}_{k};{k\in\mathbb{Z}%
}\}$ be an i.i.d. sequence of random variables with $\mathbb{P}({\varepsilon
}_{1}=0)=\mathbb{P}({\varepsilon}_{1}=1)=1/2$ and let
\[
Y_{n}=\sum_{k=0}^{\infty}2^{-k-1}{\varepsilon}_{n-k}\quad and\quad
X_{n}=g(Y_{n})-\int_{0}^{1}g(x)dx\,,
\]
where $g\in\mathbb{L}_{2}(0,1)$, $(0,1)$ being equipped with the Lebesgue measure.

By applying Proposition 3 in Maxwell and Woodroofe (2000) for verifying
condition (\ref{MW}), we see that if $g\in\mathbb{L}_{2}(0,1)$ satisfies
\begin{equation}
\int_{0}^{1}\int_{0}^{1}[g(x)-g(y)]^{2}\frac{1}{|x-y|}(\log[\log\frac
{1}{|x-y|}])^{t}dxdy<\infty\label{double}%
\end{equation}
for some $t>1$, then (\ref{MW}) is satisfied and therefore (\ref{asCLT}) and
(\ref{asCIP}) hold with $\sigma^{2}=\lim_{n\rightarrow\infty}$var$(S_{n})/n$.
A concrete example of a map satisfying (\ref{double}), pointed out in
Merlev\`{e}de-Peligrad-Utev, 2006 is%
\[
g(x)=\frac{1}{\sqrt{x}}\frac{1}{[1+\log(2/x)]^{4}}\sin(\frac{1}{x}%
)\;,\;0<x<1.
\]

\subsection{Random walks on orbits of probability preserving transformation}

\quad The following example was considered in Derriennic-Lin (2007) and also
in Cuny-Peligrad (2012). Let us recall the construction.

Let $\tau$ be an invertible ergodic measure preserving transformation on
$(S,{\mathcal{A}},\pi)$, and denote by $U$, the unitary operator induced by
$\tau$ on $\mathbb{L}_{2}(\pi)$. Given a probability $\nu=(p_{k}%
)_{k\in\mathbb{Z}}$ on $\mathbb{Z}$, we consider the Markov operator $Q$ with
invariant measure $\pi$, defined by
\[
Qf=\sum_{k\in\mathbb{Z}}p_{k}f\circ\tau^{k},\quad\text{for every }%
f\in\mathbb{L}_{1}(\pi).
\]

This operator is associated to the transition probability
\[
Q(x,A)=\sum_{k\in\mathbb{Z}}p_{k}\mathbf{1}_{A}(\tau^{k}s),\qquad s\in
S,A\in{\mathcal{A}}.
\]

We assume that $\nu$ is ergodic, i.e. the the group generated by
$\{k\in\mathbb{Z}:p_{k}>0\}$ is $\mathbb{Z}$. As shown by Derriennic-Lin
(2007), since $\tau$ is ergodic, $Q$ is ergodic too. We assume $\nu$ is
symmetric implying that the operator $Q$ is symmetric.

Denote by $\Gamma$ the unit circle. Define the Fourier transform of $\nu$ by
$\varphi(\lambda)=\sum_{k\in\mathbb{Z}}p_{k}\lambda^{k}$, for every
$\lambda\in\Gamma$. Since $\nu$ is symmetric, $\varphi(\lambda)\in
\lbrack-1,1]$, and if $\mu_{f}$ denotes the spectral measure (on $\Gamma$) of
$f\in\mathbb{L}_{2}(\pi)$, relative to the unitary operator $U$, then, the
spectral measure $\rho_{f}$ (on $[-1,1]$) of $f$, relative to the symmetric
operator $Q$ is given by%
\[
\int_{-1}^{1}\psi(s)\rho_{f}(ds)=\int_{\Gamma}\psi(\varphi(\lambda))\mu
_{f}(d\lambda),
\]
for every positive Borel function $\psi$ on $[-1,1]$. Condition (\ref{log}) is
verified under the assumption
\[
\int_{\Gamma}\frac{(\log^{+}|\log(1-\varphi(\lambda))|)^{2}}{1-\varphi
(\lambda)}\mu_{f}(d\lambda)<\infty.
\]
and therefore (\ref{PxCLT}) holds.

When $\nu$ is centered and admits a moment of order $2$ (i.e. $\sum
_{k\in\mathbb{Z}}k^{2}p_{k}<\infty$), Derriennic and Lin (2007) proved that
the condition $\int_{\Gamma}\frac{1}{|1-\varphi(\lambda)|}\mu_{f}%
(d\lambda)<\infty$, is sufficient for (\ref{PxCLT}).

\medskip

Let $a\in\mathbb{R}-\mathbb{Q}$, and let $\tau$ be the rotation by $a$ on
$\mathbb{R}/\mathbb{Z}$. Define a measure $\sigma$ on $\mathbb{R}/\mathbb{Z}$
by $\sigma=\sum_{k\in\mathbb{Z}}p_{k}\delta_{ka}$. For that $\tau$, the
canonical Markov chain associated to $Q$ is the random walk on $\mathbb{R}%
/\mathbb{Z}$ of law $\sigma$. In this setting, if $(c_{n}(f))$ denotes the
Fourier coefficients of a function $f\in\mathbb{L}_{2}(\mathbb{R}/\mathbb{Z}%
)$, condition (\ref{log}) reads
\[
\sum_{n\in\mathbb{Z}}\frac{(\log^{+}|\log(1-\varphi(\mathrm{e}^{2i\pi
na}))|)^{2}|c_{n}(f)|^{2}}{1-\varphi(\mathrm{e}^{2i\pi na})}<\infty.
\]

\subsection{CLT started from a point for a Metropolis Hastings algorithm.}

In this subsection we mention a standardized example of a stationary
irreducible and aperiodic Metropolis-Hastings algorithm with uniform marginal
distribution. This type of Markov chain is interesting since it can easily be
transformed into Markov chains with different marginal distributions. Markov
chains of this type are often studied in the literature from different points
of view. See, for instance Rio (2009).

Let $E=[-1,1]$ and let $\upsilon$ be a symmetric atomless law on $E$. The
transition probabilities are defined by
\[
Q(x,A)=(1-|x|)\delta_{x}(A)+|x|\upsilon(A),
\]
where $\delta_{x}$ denotes the Dirac measure. Assume that $\theta=\int
_{E}|x|^{-1}\upsilon(dx)<\infty$. Then there is a unique invariant measure
\[
\pi(dx)=\theta^{-1}|x|^{-1}\upsilon\,(dx)
\]
and the stationary Markov chain $(\gamma_{k})$ generated by $Q(x,A)$ and $\pi$
is reversible and positively recurrent, therefore ergodic.

\begin{thm}
Let $f$ be a function in $\mathbb{L}_{2}^{0}(\pi)$ satisfying $f(-x)=-f(x)$
for any $x\in E.$ Assume that for some positive $t$, $|f|\leq g$ on
$[-t,t]\ $where $g$ is an even positive function on $E$ such that $g$ is
nondecreasing on $[0,1]$, $x^{-1}g(x)$ is nonincreasing on $[0,1]$ and
\begin{equation}
\int_{\lbrack0,1]}[x^{-1}g(x)]^{2}dx<\infty. \label{integral}%
\end{equation}
Define $X_{k}=f(\gamma_{k}).$ Then (\ref{PxCLT}) holds.
\end{thm}

\textbf{Proof}. Because the chain is Harris recurrent if the annealed CLT
holds, then the CLT also holds for any initial distribution (see Chen, 1999),
in particular started at a point. Therefore it is enough to verify condition
(\ref{CLT-rev}). Denote, as before, by $\mathbb{E}^{x}$ the expected value for
the process started from $x\in E.$ We mention first relation (4.6) in Rio
(2009). For any $n\geq1/t$%
\[
|\mathbb{E}^{x}(S_{n}(g))|\leq ng(1/n)+t^{-1}|f(x)|\text{ for any }x\in
\lbrack-1,1].
\]
Then%
\[
|\mathbb{E}^{x}(S_{n}(g))|^{2}\leq2[ng(1/n)]^{2}+2t^{-2}|f(x)|^{2}\text{ for
any }x\in\lbrack-1,1],
\]
and so, for any $n\geq1/t$%
\[
||\mathbb{E}^{x}(S_{n})||_{2,\pi}^{2}\leq2[ng(1/n)]^{2}+2t^{-2}||f(x)||_{2,\pi
}^{2}.
\]
Now we impose condition (\ref{CLT-rev}) involving $||\mathbb{E}^{x}%
(S_{n})||_{2,\pi}^{2},$ and note that
\[
\sum\nolimits_{n}\frac{[ng(1/n)]^{2}}{n^{2}}<\infty\text{ if and only if
(\ref{integral}) holds.}%
\]

\bigskip

\textbf{Acknowledgement.} The author would like to thank the referees for
carefully reading the manuscript and for many useful suggestions that improved
the presentation of this paper.

\end{document}